\documentclass[a4paper,leqno,11pt,final]{amsproc}
\usepackage{amscd,amssymb,latexsym,amsxtra,amsfonts,amsthm,amsmath,verbatim}

\newcommand{\Q}{\mathbb Q}

\newcommand{\A}{\mathbb A}

\newcommand{\K}{\mathcal K}

\newcommand{\fA}{\mathfrak a}
\newcommand{\fS}{\mathfrak S}
\newcommand{\fO}{\mathcal O}
\newcommand{\Sym}{\rm Sym}

\makeatletter
\def\1@section{\@tocline{1}{4pt}{1pc}{}{}}
\def\1@subsection{\@tocline{2}{0pt}{2pc}{5pc}{}}
\makeatother

\begin{document}
\title{\bf Quotients of $E^n$ by $\fA_{n+1}$ and Calabi-Yau manifolds}
\author{Kapil Paranjape}
\address{The Institute of Mathematical Sciences, C.I.T. Campus,
Chennai \linebreak 600 113, India.} \email{kapil@imsc.res.in}
\author{Dinakar Ramakrishnan}
\address{Department of Mathematics, California Institute of
Technology, \linebreak Pasadena, CA 91125, USA.}
\email{dinakar@its.caltech.edu} 
\begin{abstract}We give a simple construction, for $n \geq 2$, of an $n$-dimensional
Calabi-Yau variety of Kummer type by studying the quotient $Y$ of
an $n$-fold self-product of an elliptic curve $E$ by a natural
action of the alternating group $\fA_{n+1}$ (in $n+1$ variables).
The vanishing of $H^m(Y, \fO_Y)$ for $0<m<n$ follows from the lack
of existence of (non-zero) fixed points in certain representations
of $\fA_{n+1}$. For $n \leq 3$ we provide an explicit (crepant)
resolution $X$ in characteristics different from $2, 3$. The key
point is that $Y$ can be realized as a double cover of $\mathbb
P^n$ branched along a hypersurface of degree $2(n+1)$.
\end{abstract}
\date{}

\maketitle

\section*{Introduction}

A {\it Calabi-Yau manifold} over a field $k$ is a smooth
projective variety $X$ of dimension $n$ such that
\begin{enumerate} \item[(CY1)] The canonical bundle $\K_X$ is
trivial; \, and \item[(CY2)] $H^m(X, \fO_X) = 0$ for all
(strictly) positive $m <n$. \end{enumerate} The condition (CY2) is
equivalent (for smooth $X$) to requiring that $h^{m,0}(X)=0$ for
all $m$ such that $0 < m < n$. Classically, a Calabi-Yau manifold
of dimension $n\geq 2$ is a {\it complex K\"ahler $n$-manifold
with finite $\pi_1$ (fundamental group) and SU$(n)$-holonomy}
([V]). The equivalence of the definitions is given by a {\it
theorem of S.T.~Yau}.

It will be necessary for us to allow $X$ to have mild
singularities. By a {\it Calabi-Yau variety}, we will mean a
projective variety $X/k$ on which the canonical bundle $\K_X$ is
defined such that the conditions (CY1), (CY2) hold. More
precisely, we will want such an $X$ to be normal and
Cohen-Macaulay, so that the dualizing sheaf $\K_X$ is defined,
with the singular locus in codimension at least $2$, so that
$\K_X$ defines a Weil divisor; finally $X$ should be
$\Q$-Gorenstein, so that a power of $\K_X$ will represent a
Cartier divisor.

Clearly, every Calabi-Yau manifold of dimension $1$ is an elliptic
curve, while in dimension $2$ it is a $K3$-surface. Abelian
varieties, which generalize the elliptic curve in one direction,
have trivial canonical bundles but they have non-trivial
$h^{m,0}(X)$ for $m < n$.

A classical construction of Kummer associates a $K3$ surface to an
abelian surface $A$ by starting with the quotient of $A$ by the
involution $\iota: x \to -x$, and then blowing up the sixteen
double points, each of which corresponds to a point of order $2$
on $A$. When $E$ is an elliptic curve with CM (short for {\it
complex multiplication}) by $\Q[\sqrt{-3}]$, there is a
construction of a Calabi-Yau $3$-fold arising as a resolution of a
quotient of $E\times E \times E$.

The object of this Note is to present a simple construction of a
Calabi-Yau variety {\it of Kummer type} by starting with an
$n$-fold product $E \times \dots \times E$ of an elliptic curve
$E$, and then taking a quotient under an action of the alternating
group $\fA_{n+1}$. For general $n$ this will lead, under a
suitable (crepant) resolution predicted by a standard conjecture,
to a Calabi-Yau manifold. We can do this unconditionally for
$n\leq 3$, where after getting to a local problem, one can appeal
to known results -- \cite{Ro}, for example. But we take a {\it
direct geometric approach} to arrive at the smooth resolution, and
this can at least partly be carried out for arbitrary $n$. This
construction works whether or not $E$ has CM, and it will be used
in a forthcoming paper (\cite{PaRa}). Already for $n=2$, it is
different from the classical Kummer construction (\cite{Hud});
what we do here is to divide $E \times E$ by the cyclic group of
automorphisms of order $3$ generated by $(x,y) \to (-x-y,x)$.
However, the realization of the $K3$ surface as the double cover
of $\mathbb P^2$ branched along the dual of a plane cubic has
arisen in the previous works of Barth, Katsura, and others.

The construction appears to work for $n=4$ and also over families
of elliptic curves. We plan to take up these matters elsewhere.

For $n=3$, our Calabi-Yau variety is realized as a double cover of
$\mathbb P^3$ branched along an {\it irreducible} octic surface.
For other examples of double constructions, with highly reducible
branch locus, see \cite{CyM} (and the references therein).

The second author would like to thank the organizers of the {\it
International Conference on Algebra and Number theory} for
inviting him to come to Hyderabad, India, during December 2003,
and participate in the conference. This Note elaborates on a small
part of the actual lecture he gave there, describing the ongoing
joint work with the first author.

We would like to thank Slawomir Cynk and Klaus Hulek for reading
an earlier version of this paper carefully and for pointing out an
error there in the resolution for $n=3$, and moreover for
suggesting an elegant way to fix it, which we have used in section
2.2. We would also like to thank Igor Dolgachev for interesting
suggestions for further work which we plan to take up in a sequel.

\bigskip

\section{The construction}

\medskip

Let $E$ be an elliptic curve over a field $k$ with identity $0$,
and $n \geq 1$ an integer. Put
$$
\tilde Y: = \, \left\{\tilde y=(y_1, \dots, y_{n+1}) \in E^{n+1} \,
\vert \, \sum_{j=1}^{n+1} y_j = 0\right\}\leqno(1.1)
$$
Clearly, we have an isomorphism
$$
\varphi: \tilde Y \, \rightarrow \, E^n,\leqno(1.2)
$$
given by $\tilde y \to (y_1, \dots , y_n)$.

Note that the action of the alternating group $\fA_{n+1}$ on
$E^{n+1}$ preserves $\tilde Y$. Put
$$
Y: = \, \tilde Y/\fA_{n+1}. \leqno(1.3)
$$

This variety is defined over $k$, but is singular. Denote by $\pi:
\tilde Y \to Y$ the quotient map and by $Z$ the singular locus in
$Y$. If we set
$$
\tilde Z: = \, \left\{\tilde y \in \tilde Y \, \vert \, \exists g
\in \fA_{n+1}, g \ne 1, \, {\rm s.t.} \, \, g\tilde y = \tilde
y\right\},\leqno(1.4)
$$
namely the set of points in $\tilde Y$ with non-trivial stabilizers
in $\fA_{n+1}$, we obtain
$$
Z \, \subset \, \pi(\tilde Z).\leqno(1.5)
$$
If $n=2$, for example, the action of $\fA_3$ on $E \times E$ (via
$\varphi$) is generated by $(x,y)\to (-x-y,x)$, which shows that the
fixed point set is $\{(x,x)\in E\times  E \, \vert \, 3x=0\}$.

\medskip

\noindent{\bf Theorem} \, \it We have the following (for $n \geq
2$):

\noindent(a) \, $Y$ is a Calabi-Yau variety i.e., $\K_Y$ is
defined with
\begin{enumerate}
\item[(i)] $\K_Y$ is trivial
\item[(ii)] $H^m(Y, \fO_Y) = 0$ for all $m$ such that $0 < m < n$
\end{enumerate}

\noindent(b) \, If $n \leq 3$ and $k$ algebraically closed of
characteristic zero or $p \nmid 6$, there exists a smooth resolution
$p: X \to Y$ such that $X$ is Calabi-Yau. \rm

\medskip

{\bf Proof of Theorem, part ${\bf (a)}$}: \, We need the following:

\medskip

\noindent{\bf Proposition A} \, \it Consider the morphism $\pi:
\tilde Y \, \rightarrow \, Y = \tilde Y/\fA_{n+1}$. Then $\pi$ is
finite, surjective and separable. Moreover, the natural homomorphism
$$
\fO_Y \, \rightarrow \, \pi_\ast(\fO_{\tilde Y})^{\fA_{n+1}}
$$
is an isomorphism. \rm

\medskip

{\bf Proof of Proposition A}. \, In view of the Theorem in chap.
II, sec. 7 of [Mu], it suffices to prove that for any point
$\tilde y = (y_1, \dots, y_{n+1})$ in $\tilde Y$, the orbit
$O(\tilde y)$ is contained in an affine open subset of $\tilde Y$.
(In fact one should properly appeal to this Theorem of Mumford to
already know that the algebraic quotient $Y$ exists and is
unique.) Now by definition, $y_{n+1}=-\sum_{j=1}^n y_j$ for any
$\tilde y = (y_1, \dots, y_{n+1})$. Pick any affine open set $U$
in $E$ which avoids the points $\{y_1, \dots, y_{n+1}\}$. Then
$U^n$ is an affine open subset of $E^n$, ad the orbit $O(\tilde
y)$ is contained in the affine open subset $\varphi^{-1}(U^n)$ of
$\tilde Y$. Done.

\qed

\medskip

Put
$$
W: = \, H^1(E, \fO_E) \simeq k\leqno(1.6)
$$
and
$$
W_{m,n} \, = \, \Lambda^m(W^{\oplus n}) \simeq H^m(E^n, \fO_{E^n}).
$$
In view of Proposition A and the isomorphism $\phi$, we are led to
look for fixed points of the action of $\fA_{n+1}$ on $W_{m,n}$.
To be precise, our Theorem will be a consequence of the following

\medskip

\noindent{\bf Proposition B} \, \it Fix $n \geq 2$. Let $k$ have
characteristic zero or $p \nmid (n!/2)$. Then for every integer
$m$ such that $0 < m < n$,
$$
W_{m,n}^{\fA_{n+1}} \, = \, 0.
$$
\rm

\medskip

{\bf Proof of Proposition B}. \,First consider the simple case
${\bf n=2}$. Here the only possibility is $m=1$. The group $\fA_3$
is generated by the $3$-cycle $\begin{pmatrix}1 & 2 &
3\end{pmatrix}$, which sends $(w_1, w_2) \in W_{1,2}$ to
$(-w_1-w_2,w_1)$ and is represented by the matrix
$\begin{pmatrix}-1 & -1\\ 1 & 0
\end{pmatrix}$. Since char$(k) \ne 3$, the eigenvalues are the two
non-trivial cube roots of unity, implying that there is no fixed
point in $W_{1,2}$.

So we may take ${\bf n \geq 3}$ and assume by induction that the
Proposition is true for $n-1$. Put
$$
W'_{1,n} \, = \, \{w=(w_1,\dots, w_{n+1}) \in W^{n+1} \, \vert \,
w_1=0, \sum_{j=2}^{n+1} w_j = 0\},\leqno(1.7)
$$
$$
L \, = \, \{w=(w_1,\dots,w_{n+1}) \in W^{n+1} \, \vert \, w_1=n,
w_j=-1 \, \, \forall j \geq 2\},
$$
and
$$
G': = \{g \in \fA_{n+1} \, \vert \, g(w_1) = w_1\}.
$$
Then there are canonical, compatible identifications $W'_{1,n}
\simeq W_{1,n-1}$ and $G' \simeq \fA_n$, and so by induction,
$$
\Lambda^j(W'_{1,n})^{G'} = 0 \quad {\rm if} \quad 0<j<n-1.
\leqno(1.8)
$$
Moreover, since $W_{1,n}$ identifies with the $\fA_{n+1}$-space of
vectors $(w_1, \dots, w_{n+1}) \in W^{n+1}$ such that $\sum_j w_j
= 0$, we get a $G'$-stable decomposition
$$
W_{1,n} \, = \, W'_{1,n} \oplus L,\leqno(1.9)
$$
with $G'$ acting trivially on the line $L$. This furnishes, by
taking exterior powers, $G'$-isomorphisms for all positive
integers $m\leq n-1$,
$$
\Lambda^m(W_{1,n}) \, \simeq \, \Lambda^m(W'_{1,n}) \, \oplus \,
\Lambda^{m-1}(W'_{1,n}) \otimes L. \leqno(1.10)
$$
We then get, by the inductive hypothesis,
$$
W_{1,n}^{G'} \, = \, 0 \quad {\rm if}\quad 1 < m <
n-1.\leqno(1.11)
$$

So it suffices to prove the Proposition for $m=1$ and $m=n-1$. We
will be done once we show the following

\medskip

\noindent{\bf Lemma 1.12} \, \it For $n \geq 3$, the
representation $\rho$ of $G = \fA_{n+1}$ on $W_{1,n}$ is
irreducible. Moreover,
$$
W_{n-1,n} \, \simeq \, W_{1,n}^\vee \otimes {\rm det}(\rho),
$$
where the superscript $\vee$ denotes taking the contragredient. \rm

\medskip

{\bf Proof of Lemma}. \, Assume for the moment the irreducibility
of $\rho$. As dim$W_{1,n}=n$, there is a natural, non-degenerate
$G$-pairing
$$
\Lambda^{n-1}(W_{1,n}) \times W_{1,n} \, \rightarrow \,
\Lambda^n(W_{1,n}), \, (\alpha, w) \to \alpha\wedge w,
$$
and $G$ acts on the one-dimensional space on the right by
det$(\rho)$, identifying the representation of $G$ on $W_{n-1,n}$
with $\rho^\vee \otimes {\rm det}(\rho)$.

\medskip

All that remains now is to check the irreducibility of $\rho$. For
this note that the action of $G=\fA_{n+1}$ on $V= k^{n+1}$ is
doubly transitive, and so by [CuR], this permutation
representation $\pi$, say, decomposes as the direct sum of the
trivial representation and an irreducible representation of $G$,
which must be equivalent to $\rho$. But here is an explicit
argument. Since $G'$ is the stabilizer of $(1,0,\dots,0)$ in $V$,
we see that $\pi$ is the representation induced by the trivial
representation of $G'$. On the other hand, the double coset space
$G'\backslash G/G'$ has exactly two elements, again implying, by
Mackey, that the complement of $1$ in $\pi$ is irreducible. Done.

\qed

\medskip

Now we turn to the question of triviality of $\K_Y$. As $\tilde Y$
is an abelian variety, $\K_{\tilde Y}$ is trivial. The quotient
$Y$ is Cohen-Macaulay, being a finite group quotient of a smooth
variety. It is normal with the singular locus in codimension $2$,
and is $\Q$-Gorenstein. $\K_Y$ identifies with the line bundle on
$Y$ defined by the $\fA_{n+1}$-invariance of $\K_{\tilde Y}$.
Moreover, there is a section of $\K_{\tilde Y}$ which is
invariant. This gives a section of $\K_Y$ over $Y$, showing the
triviality of $\K_Y$. (This argument will not work if we divide by
the full symmetric group, because then any transposition will act
by $-1$ upstairs, and the section will not be invariant.)
Alternatively, we will show below that $Y$ is a double cover of
${\mathbb P}^n$ branched along a hypersurface of degree $2n+2$,
again implying that $\K_Y$ is trivial.

\medskip

We have now proved part (a) of our Theorem.

\bigskip

\section{Resolution}

\medskip

Now we will show how to deduce part (b) of Theorem. To begin,
since the variety $Y$ constructed above is an orbifold, a standard
conjecture predicts that there will be a smooth resolution
$$
p: X \rightarrow Y
$$
which is {\it crepant}, i.e., that the canonical bundle of $X$ has
for image the canonical bundle of $Y$ (under $p_\ast$) and is thus
trivial. For $n \leq 3$ this can be achieved by making use of
[Ro], but we will take a different tack.

\medskip

Now consider the natural action of the symmetric group $\fS_{n+1}$
on $E^{n+1}$, the product of $n+1$ copies of $E$. The addition map
$E^{n+1}\to E$ is stable under the action of $\fS_{n+1}$ and thus
we obtain a map $\Sym^{n+1}(E)\to E$, where the former denotes the
quotient of $E^{n+1}$ by the action.

The space $\Sym^{n+1}(E)$ can also be identified with the space of
effective divisors of degree $n+1$ on $E$ and under this
identification the above map can be understood as follows. Let $o$
denote the origin in $E$. For each point $p$ in $E$ the fibre of
the map consists of all divisors in the linear system
$|n[o]+[p]|$. In particular, when $p=o$ we see that the fibre
consists of all divisors in the linear system $|(n+1)[o]|$.

From the point of view of quotients the fibre over $o$ is the
quotient by the action of $\fS_n$ of the space
\[ \tilde Y \, = \,  \{ (p_0,\dots,p_n) | p_0 + \cdots + p_n = 0 \} \]
We are interested in the quotient $Y$ of this space by the
alternating group $\fA_{n+1}$. Thus, $Y$ can be expressed as a
double cover of the linear system $|(n+1)[o]|$ branched along the
locus of divisors of the form $2[p]+[p_2]+\cdots+[p_n]$.

When $n$ is at least 2 the linear system $|(n+1)[o]|$ gives an
embedding of $E$ into the dual projective space $|(n+1)[o]|^{*}$.
The locus of special divisors as considered above is then
identified with the dual variety of $E$; i.~e.,\ the variety that
consists of all hyperplanes that are tangent to $E$. It is well
known that this dual variety has degree $2(n+1)$, which follows
for example from the Hurwitz genus formula giving the number of
ramification points for a map $E\to{\mathbb P}^1$ of degree $n+1$.

Since $Y$ is a double cover of ${\mathbb P}^n$ branched along this
hypersurface of degree $2(n+1)$, as claimed above, implying the
triviality of $\K_Y$. In order to find a good resolution of $Y$ it
is sufficient to understand the singularities of the dual variety.

\subsection{The case $n=2$}
Here we have the dual of the familiar embedding of $E$ as a cubic
curve in ${\mathbb P}^2$. This curve has 9 points of inflection
and no other unusual tangents. It follows from the usual theory
that the dual curve is a curve with 9 cusps and no other
singularities. Thus $Y$ is the double cover of ${\mathbb P}^2$
branched along such a curve. To resolve $Y$ it is enough to
resolve over each cusp individually.

Thus we consider the simpler case of resolving the double cover of
$W\to\A^2$ branched along the curve defined by $y^2-x^3$; the
variety $W$ is a closed subvariety of $\A^3$ defined by
$z^2-y^2+x^3$ with the projection to the $(x,y)$ plane providing the
double covering. One checks easily that the blow-up of the maximal
ideal $(x,y,z)$ gives a resolution of singularities. Moreover, this
blow-up is a double cover of the blow-up of $\A^2$ at the maximal
ideal $(x,y)$. Since the exceptional divisor in the first case is a
$(-1)$ curve, it follows that the exceptional divisor in the blow-up
of $W$ is a $(-2)$ curve.

Let $X\to Y$ be the result of blowing-up the nine singular points
in $Y$ that lie over the cusps of the dual curve; as seen above
$X$ is smooth. From the adjunction formula we see that $\K_X$
restricts to the trivial divisor on each exceptional divisor;
hence $\K_X$ is the pull-back of the $\K_Y$. The usual theory of
double covers shows us that $\K_Y$ is trivial and $Y$ is
simply-connected. Thus the same properties hold for $X$ as well.
In other words we have shown that $X$ is a K3 surface.

\subsection{The case $n=3$}
In this case $E$ is embedded as the complete intersection of a
pencil of quadrics in $\mathbb P^3$. Recall that we have assumed
that the characteristic does not divide $6$.

A point of the dual variety $D$ corresponds to a plane that
contains a tangent line. Thus each point on $E$ determines a
pencil of such points. Equivalently, if $P\subset E\times
(P^3)^{*}$ denotes the projective bundle on $E$ that consists of
pairs $(p,\pi)$ where $\pi$ is a plane in $P^3$ that is tangent to
$E$ at $p$, $D$ is the image of $P$ under the natural projection
to $(P^3)^*$ which is a surface of degree 8. For notational
convenience let the origin of the group law on $E$ be chosen to be
a point $o$ such that the linear system is $4[o]$. The fibre of
$P$ over a point $p$ can then also be described as the collection
of all divisors $D=[q]+[r]$ of degree two such that $2p+q+r=o$ in
the group law.

Let $a$ be a point of order two in $E$. Then for each point $p$ in
$E$ we can consider the point $2[a-p]$ in the fibre of $P$ over
$p$; this gives a section $\sigma_{a}:E\to P$ and we denote the
image in $P$ as $E_a$. This gives us four disjoint curves in $P$.
Under the composite map $E_a\to P\to D$, the point $\sigma_a(p)$
and $\sigma_a(a-p)$ are both sent to the hyperplane that
intersects $E$ in $2[p]+2[a-p]$, so the image $C_a$ of $E_a$ in
$D$ is the quotient of $E_a$ by the involution $p \mapsto a-p$;
thus $C_a$ is isomorphic to $P^1$. Moreover, $D$ has a transverse
ordinary double point along the general point of $C_a$.

Let $p$ be any point of $E$. Then $[-3p]+[p]$ is a point in the
fibre of $P$ over $E$; this gives a section $\tau:E\to P$ and we
denote the image in $P$ as $F$. The composite map $F\to P\to D$ is
a one-to-one since the hyperplane section of the type $3[p]+[q]$
uniquely determines the point $p$; let $G$ denote the image of $F$
in $D$. One notes that $D$ has a transverse cusp along the general
point of $G$.

Let $b$ be a point of order 4 on $E$ and consider the point $a=2b$
which is a point of order 2 on $E$. We see that
\[ \tau(b) = [-3b]+[b]=2[b] = 2[a-b] = \sigma_a(b) \]
The curve $E_a$ thus intersects the curve $F$ in the four points
of order 4 which are "half" of $a$. As $a$ varies we obtain 16
points on $P$ lying over the 16 points of order 4 on $E$. The
singularity of $D$ at the image of these sixteen points can be
described as follows in suitable local coordinates $x$, $y$ and
$z$. The curve $C_a$ is described by $y=z+h=0$ and the curve $G$
is described by $x^3+y^2=z+h'=0$, where $h$ and $h'$ consists of
terms of degree at least two and are {\em distinct} (this is
important in the next paragraph). Moreover, the Jacobian ideal of
$D$ is the intersection of the ideals $(y,z+h)$ and
$(x^3+y^2,z+h')$ that define these two curves.

Let $f:U\to P^3$ be the the smooth threefold obtained by blowing
up along the curves $C_a$. Let $A_a$ denote the exceptional locus
in $U$ over the curve $C_a$. The canonical divisor of $U$ is
$f^*\fO(-4)\otimes \fO(\sum A_a)$. Since $D$ has an ordinary
double point along $C_a$ we see that the strict transform $D'$ of
$D$ is linearly equivalent to $f^*\fO(8)\otimes \fO(-2 \sum A_a)$.
Moreover, $D'$ is only singular along the strict transform $G'$ of
$G$. Finally, from the above local description it follows that the
surface $z+h'=0$ which $G'$ lies on {\em is} blown up at the
origin $(x,y)$ under $f$. It follows that $G'$ is smooth and $D'$
has a transverse cusp along it.

Let $g:Z\to U$ be the smooth threefold obtained by  blowing up
along $G'$. Let $B$ denote the exceptional locus in $Z$ over $G$
and (by abuse of notation) let $A_a$ again denote the strict
transform of the divisors $A_a$ in $U$. The canonical divisor of
$Z$ is $g^*f^*\fO(-4)\otimes \fO(B+\sum A_a)$. Since $D'$ is a
singularity of multiplicity two along $G'$ we see again that its
strict transform differs from its total transform by $2B$; thus
the strict transform $D''$ of $D'$ is linearly equivalent to
$g^*f^*\fO(8)\otimes \fO(2B-2\sum A_a)$. Finally, we see that
$D''$ is smooth as well. Thus the double cover $Y\to Z$ along
$D''$ is a smooth threefold with trivial canonical bundle.

\qed

\vskip 0.3in

\end{document}